\documentclass[MSNbibl,number,citesort,seceqn,dvips]{arxbj}

\aid{0}
\volume{19}
\issue{4}
\pubyear{2013}
\firstpage{1350}
\lastpage{1377}
\doi{10.3150/12-BEJSP01}

\makeatletter

\newtheorem{thmm}{Theorem}[section]
\newtheorem{pro}{Proposition}[section]
\newtheorem{cor}{Corollary}[section]
\newremark{rem}{Remark}[section]
\newproclaim{Def}{Definition}[section]
\newcommand{\U}{\underline{U}}
\newcommand{\R}{{\mathbb R}}
\newcommand{\N}{{\mathbb N}}
\newcommand{\e}{\varepsilon}
\renewcommand{\l}{\lambda}
\renewcommand{\a}{\alpha}
\newcommand{\li}[2]{\lim_{#1\rightarrow #2}}
\newcommand{\I}[4]{\int_{#1}#2 #3(\mathrm{d}#4)}
\newcommand{\Io}[4]{\int_{#1}^{#2}#3 \,\mathrm{d}#4}
\newcommand{\PR}[2]{(#1,#2)}
\newcommand{\st}{\stackrel{\mathrm{(law)}}{=}}
\makeatother

\begin{document}
\begin{frontmatter}

\title{On the Mellin transforms of the perpetuity and the remainder
variables associated to a~subordinator}
\runtitle{Mellin transform of perpetuity variables}

\begin{aug}
\author[1]{\fnms{Francis} \snm{Hirsch}\corref{}\thanksref{1}\ead[label=e1]{francis.hirsch@univ-evry.fr}}%
\and
\author[2]{\fnms{Marc} \snm{Yor}\thanksref{2,3}\ead[label=e2]{deaproba@proba.jussieu.fr}}
\runauthor{F. Hirsch and M.Yor} 
\address[1]{Laboratoire d'Analyse et Probabilit\'{e}s,
Universit\'{e} d'\'{E}vry -- Val d'Essonne, Boulevard F. Mitterrand,
\mbox{F-91025} \'{E}vry Cedex,
France.
\printead{e1}}
\address[2]{Laboratoire de Probabilit\'{e}s et Mod\`{e}les
Al\'{e}atoires,
Universit\'{e} Paris VI et VII, 4 Place Jussieu~-- Case 188,
\mbox{F-75252} Paris Cedex 05,
France.
\printead{e2}}
\address[3]{Institut Universitaire de France}

\end{aug}


%
\begin{abstract}
Results about the laws of the perpetuity and remainder variables
associated to a subordinator are presented, with particular emphasis on
their Mellin transforms, and multiplicative infinite divisibility
property. Previous results by Bertoin--Yor (\textit{Electron. Commun. Probab.}
\textbf{6} (2001) 95--106) are incorporated in
our discussion; important examples when the subordinator is the inverse
local time of a diffusion are exhibited. Results of Urbanik
(\textit{Probab. Math. Statist.} \textbf{15} (1995) 493--513) are
also discussed in detail; they appear to be too little known, despite
the fact that quite a few of them have priority upon other works in
this area.
\end{abstract}

%

\end{frontmatter}

\section{Introduction}\label{s1}
Let $\PR{\xi_l}{l\geq0}$ denote a subordinator, with Laplace--Bernstein
exponent $\Phi\not\equiv0$:
\[
\mathbb{E}\bigl[\exp(-s \xi_l)\bigr]=\exp\bigl(-l \Phi(s)\bigr).
\]
We take the slightly unusual notation $l$ for the time parameter,
instead of $t$, because we have in mind, among other examples, those
$\PR{\xi_l}{l\geq0}$ which are inverse local times, i.e.:
\[
\xi_l=\inf\{t;L_t>l\},
\]
where $\PR{L_t}{t\geq0}$ is (a choice of) the local time at 0 for some
1-dimensional diffusion (see, in particular, Subsection \ref{ss44}, for
the case where the diffusion is a \emph{radial Ornstein--Uhlenbeck
process}).


The present paper is a general study in the spirit of Bertoin--Yor \cite
{BY}. So,
we first recall the main results of Bertoin--Yor \cite{BY} (see also,
for property i), Carmona et al. (\cite{CPY}, Proposition~3.3)):
\begin{enumerate}[iii)]
\item[i)]
The \emph{perpetuity}:
\[
\mathcal{I}\equiv\mathcal{I}_{\xi}\equiv\Io{0} {\infty} {\exp (-
\xi_l)} {l}
\]
has integral moments of all orders, which determine its law:
%
\begin{equation}
\label{e11} \mathbb{E}\bigl[\mathcal{I}^n\bigr]=\frac{n!}{\Phi(1)\cdots\Phi(n)},\ n
\geq1.
\end{equation}
\item[ii)] There exists a random variable $\mathcal{R}\equiv\mathcal
{R}_{\xi
}$ (the \emph{remainder}) whose law is also determined by its integral moments:
%
\begin{equation}
\label{e12} \mathbb{E}\bigl[\mathcal{R}^n\bigr]=\Phi(1)\cdots
\Phi(n),\ n\geq1.
\end{equation}
In Bertoin--Yor \cite{BY}, the variable $\mathcal{R}$ is defined via the
relation:
\[
\mathbb{E} \bigl[\exp(-t \mathcal{R}) \bigr]=\mathbb{E} \biggl[\frac
{1}{X_t}
\biggr],\ t\geq0
\]
where $X$ denotes the Lamperti process, starting from 1, associated
with the subordinator~$\xi$.
\item[iii)] There is the factorization:
%
\begin{equation}
\label{e13} \mathrm{\mathbf{e}}\stackrel{\mathrm{(law)}} {=}\mathcal{I}\cdot
\mathcal{R},
\end{equation}
where $\mathrm{\mathbf{e}}$ denotes a standard exponential variable,
and, on the RHS of (\ref{e13}),
$\mathcal{I}$ and $\mathcal{R}$ are independent.
\end{enumerate}

We recall that Berg--Dur\'{a}n (\cite{BD}, Theorem 1.3) deduced that the
right-hand sides of (\ref{e11}) and~(\ref{e12}) are Stieltjes moment
sequences, from their general result:

if $\PR{f(s)}{s>0}$ is a completely monotone function, then:
\[
\frac{1}{f(1)\cdots f(n)},\ n=1,2,\ldots,\quad \mbox{is a moment-sequence}
\]
(for (\ref{e11}), take $f_1(s)=\Phi(s)/s$, and, for (\ref{e12}), take
$f_2(s)=1/\Phi(s)$).

In the sequel, we shall often refer to the following \emph{trivial example}:
\[
\Phi(s)=s;\ \xi_l=l;\ \mathcal{I}=1;\ \mathcal{R}\st{\mathbf
{e}};\ \mathbb{E}\bigl[\mathcal{R}^n\bigr]=n!,\ n\in\N.
\]

In the present paper, we discuss a number of precisions and
improvements about the previous results i) and ii):

In Section \ref{s2}, we introduce, for $r>0$, the Mellin transforms:
\[
I(r)=\mathbb{E}\bigl[\mathcal{I}^{r-1}\bigr],\ R(r)=\mathbb{E}\bigl[
\mathcal{R}^{r-1}\bigr].\vadjust{\goodbreak}
\]
Note that, in the trivial example, $R(r)=\Gamma(r),r>0$. We prove that,
in the general case, the following functional equations hold:
\[
I(r+1)=\frac{r}{\Phi(r)} I(r),\ R(r+1)=\Phi(r) R(r).
\]
This leads us to a study of the Mellin transforms $I$ and $R$,
following the same method as Artin~\cite{A} in his study of the Gamma function.

In Section \ref{s3}, we show the infinite divisibility of the random
variable $\log(\mathcal{R})$. This leads to some integral
representations of $\log(R(r))$ and $\log(I(r))$. (Note that the
integral representation of $\log(R(r)$) was also obtained by Berg
(\cite{B2}, Theorem 2.2).) As a consequence, we obtain a characterization of
the infinite divisibility and of the self-decomposability of the random
variable $\log(\mathcal{I})$, in terms of the measure whose Laplace
transform is $\Phi'/\Phi$.

In Section \ref{s4}, we discuss many examples, in particular when $\PR
{\xi_l}{l\geq0}$ is the inverse local time of a radial
Ornstein--Uhlenbeck process.

In Section \ref{s5}, we present the main results of Urbanik \cite{U},
a number of which predate some of the results found in the previous
sections.

To conclude this introduction, we note that the topic of this paper
has attracted a lot of attention since at least the mid-nineties, which
stems among other origins from applications to models in insurance,
telecommunications, and many other applied domains.

Although the papers \cite{B3}, \cite{B2}, \cite{BY}, \cite{BY2}, \cite
{CPY}, \cite{U} deal with topics studied in the present paper, we
believe the subject still warrants some new exposition. In particular,
the wealth of results found in Urbanik \cite{U} has not been properly
appreciated, perhaps because of some unusual (for probabilists)
notation adopted by Urbanik in his work.

We also note that the perpetuities involved in our discussion are a
very particular class among the family of exponential functionals of
the form:
\[
\int_0^{\infty}\exp(-\xi_{l_-}) \,\mathrm{d}
\eta_l,
\]
where the pair $(\xi,\eta)$ is a 2-dimensional L\'evy process, for
which the integral makes sense (see, e.g., Carmona et al. \cite{CPY}).

\section{Functional equations}\label{s2}
\subsection{Notation}\label{ss21} First, we fix the notation. In the sequel,
we consider a subordinator $\PR{\xi_l}{l\geq0}$ with \emph{Laplace--Bernstein exponent} $\Phi\not\equiv0$, i.e.:
\[
\forall l\geq0,\ \forall s\geq0,\quad\mathbb{E}\bigl[\exp(-s \xi_l)
\bigr]=\exp \bigl(-l \Phi(s)\bigr).
\]
One has:
%
\begin{equation}
\label{e20} \forall s\geq0,\quad\Phi(s)=a s+\I{(0,+\infty)} {\bigl(1-\mathrm
{e}^{-sx}\bigr)} {\l} {x}
\end{equation}
for some $a\geq0$ and a measure $\l$ which satisfies:
\[
\I{(0,+\infty)} {(x\wedge1)} {\l} {x}<\infty.
\]

\begin{Def}\label{d21}
In the sequel, we call \emph{Bernstein function} any such function $\Phi
$ satisfying (\ref{e20}), with $a\geq0$ and $\l$ a measure on
$(0,+\infty)$ such that $ {\I{(0,+\infty)}{(x\wedge1)}{\l
}{x}<\infty}$. In particular, in this paper, \emph{a Bernstein function
is always assumed to be null at 0}. The measure $\l$ is called the
\emph{L\'evy measure} of $\xi$ or $\Phi$.
\end{Def}
We denote by $\mu_l$ the law of $\xi_l$, and by $\rho$ the potential
measure defined by:
\[
\rho=\Io{0} {\infty} {\mu_l} {l}.
\]
We also set:
\[
\widehat{\l}(\mathrm{d}x)=x \l(\mathrm{d}x)\quad\mbox{and}\quad \overline{\l}(
\mathrm{d}x)=\l\bigl((x,+\infty)\bigr) \,\mathrm{d}x.
\]
Finally, we set: $\kappa= (a \e+\widehat{\l})\ast\rho$, where $\e$
denotes the Dirac measure at 0 and $\ast$ the convolution of measures.

The following proposition is easily proven.
%
\begin{pro}\label{p21}
The functions: $ {\Phi'(s),\frac{\Phi(s)}{s},\frac{1}{\Phi
(s)},\frac{\Phi'(s)}{\Phi(s)}}$, are completely monotone. They are
respectively the Laplace transforms of: $ a \e+\widehat{\l},a \e
+\overline{\l},\rho,\kappa$.
\end{pro}
In particular, the measure $\kappa$, which plays an important role in
the sequel, may also be defined by:
\[
\forall s>0 ,\quad\frac{\Phi'(s)}{\Phi(s)}=\I{\R_+} {\mathrm {e}^{-sx}} {
\kappa} {x}.
\]
Since $ {\li{s}{\infty}\frac{\Phi'(s)}{\Phi(s)}=0}$, one
has: $\kappa(\{0\})=0$.

The perpetuity and remainder variables are defined in Section \ref{s1}
and denoted by $\mathcal{I}$ and~$\mathcal{R,}$ respectively. We now
introduce the Mellin transforms:
\[
I(r):=\mathbb{E}\bigl[\mathcal{I}^{r-1}\bigr],\quad R(r):=\mathbb
{E}\bigl[\mathcal {R}^{r-1}\bigr],\quad \mbox{for }r>0.
\]
%
\subsection{Functional equations}\label{ss22}
%
\begin{pro}\label{p22}
The following functional equations hold:
%
\begin{eqnarray}
\label{e21} \forall r>0 ,\quad I(r+1)&=&\frac{r}{\Phi(r)} I(r),
\\
\label{e22}\forall r>0 ,\quad R(r+1)&=&\Phi(r) R(r).
\end{eqnarray}
\end{pro}
We note that equation (\ref{e21}) is a particular case of Proposition
3.1 in Carmona et~al.~\cite{CPY}.
\begin{pf}
\begin{enumerate}[1)]
\item[1)] We set, for $t\geq0$,
\[
\mathcal{I}_t=\Io{t} {\infty} {\exp(-\xi_l)} {l}= \biggl(
\Io {0} {\infty} {\exp \bigl(-(\xi_{l+t}-\xi_t)\bigr)} {l}
\biggr) \exp(-\xi_t).
\]
Hence, for $r>0$,
\begin{eqnarray*}
\frac{\mathrm{d}}{\mathrm{d}t} \mathbb{E}\bigl[\mathcal{I}_t^r
\bigr]&=&-r \mathbb{E}\bigl[\mathcal{I}_t^{r-1} \exp(-
\xi_t)\bigr]
\\
&=& - r \mathbb{E} \biggl[ \biggl( \Io{0} {\infty} {\exp\bigl(-(
\xi_{l+t}-\xi_t)\bigr)} {l} \biggr)^{r-1} \exp(-r
\xi_t) \biggr]
\\
&=&-r \mathbb{E}\bigl[\mathcal{I}^{r-1}\bigr] \exp\bigl(-t \Phi(r)
\bigr).
\end{eqnarray*}
Integrating between $0$ and $\infty$, we obtain (\ref{e21}).
\item[2)] By (\ref{e13}),
%
\begin{equation}
\label{e23} \forall r>0 ,\quad I(r) R(r)=\Gamma(r).
\end{equation}
Then, (\ref{e22}) follows easily from (\ref{e21}) and
(\ref{e23}).\quad\qed
\end{enumerate}
\noqed\end{pf}
In the trivial example, $\Phi(r)=r$, $R(r)=\Gamma(r)$, and (\ref{e22})
is the classical equation satisfied by the Gamma function. In the next
subsection, we use Artin's method based on logarithmic convexity, to
characterize the Mellin transforms $I$ and $R$.

\subsection{Logarithmic convexity}\label{ss23}
The following theorem extends the classical Bohr--Mollerup theorem (see,
for instance, Andrews et al. (\cite{AAR}, Theorem~1.9.3)).
%
\begin{thmm}\label{t21}
The Mellin transform: $r>0\longrightarrow R(r)$, is the unique function
$f$ from $(0,+\infty)$ into $(0,+\infty)$ such that:
\begin{enumerate}[iii)]
\item[i)]$f(1)=1$,
\item[ii)]$\forall r>0$, $f(r+1)=\Phi(r) f(r)$,
\item[iii)]$ f$ is log-convex (i.e., $\log f$ is convex) on $(0,+\infty)$.
\end{enumerate}
\end{thmm}
\begin{pf}
\begin{enumerate}[1)]
\item[1)] Clearly, $f=R$ satisfies properties i) and iii) by the
definition, and satisfies property ii) by Proposition \ref{p22}.
\item[2)] Let $f$ satisfying properties i), ii) and iii), $r\in(0,1]$ and
$n\geq2$. By iteration of ii), we get:
\begin{eqnarray*}
&\displaystyle f(n+r)=\Phi(r) \Phi(r+1) \cdots\Phi(r+n-1) f(r),&
\\
&\displaystyle f(n)=\Phi(1) \cdots\Phi(n-1).&
\end{eqnarray*}
Moreover, by property iii),
\[
f(n) \biggl(\frac{f(n)}{f(n-1)} \biggr)^r\leq f(n+r)\leq f(n) \biggl(
\frac{f(n+1)}{f(n)} \biggr)^r.
\]
Hence,
\[
\frac{1}{\Phi(r)} \prod_{j=1}^{n-1}
\frac{\Phi(j)}{\Phi(j+r)} \bigl(\Phi(n-1) \bigr)^r\leq f(r)\leq
\frac{1}{\Phi(r)} \prod_{j=1}^{n-1}
\frac{\Phi(j)}{\Phi(j+r)} \bigl(\Phi(n) \bigr)^r.
\]
Since $R$ also satisfies properties i), ii) and iii), we obtain:
\[
\biggl(\frac{\Phi(n-1)}{\Phi(n)} \biggr)^r\leq\frac
{f(r)}{R(r)}\leq
\biggl(\frac{\Phi(n)}{\Phi(n-1)} \biggr)^r.
\]
Letting $n$ tend to $\infty$, we get:
\[
\forall0<r\leq1,\quad f(r)=R(r).
\]
Using again property ii), we obtain by induction that $f=R$.\quad\qed
\end{enumerate}
\noqed\end{pf}
The same proof, replacing $\Phi(s)$ by $s/\Phi(s)$, yields the
following theorem.
%
\begin{thmm}\label{t22}
The Mellin transform: $r>0\longrightarrow I(r)$, is the unique function
$g$ from $(0,+\infty)$ into $(0,+\infty)$ such that:
\begin{enumerate}[iii)]
\item[i)]$g(1)=1$,
\item[ii)]$\forall r>0$, $ {g(r+1)=\frac{r}{\Phi(r)} g(r)}$,
\item[iii)]$ g$ is log-convex on $(0,+\infty)$.
\end{enumerate}
\end{thmm}
%
\subsection{Representation as limits}\label{ss24}
The next theorem generalizes the following classical formula:
%
\begin{equation}
\label{e24} \forall r>0,\quad\Gamma(r)=\li{n} {\infty}\frac{n!}{r (r+1) \cdots
(r+n-1)}
n^{r-1}.
\end{equation}

\begin{thmm}\label{t23} There is the asymptotic formula:
\[
\forall r>0 ,\quad R(r)=\li{n} {\infty} \prod_{j=0}^{n-1}
\frac{\Phi
(j+1)}{\Phi(j+r)} \bigl(\Phi(n) \bigr)^{r-1}.\vadjust{\goodbreak}
\]
\end{thmm}
\begin{pf}
\begin{enumerate}[1)]
\item[1)] We set, for $n\geq1$ and $r>0$,
\[
h(n,r)= \prod_{j=0}^{n-1}\frac{\Phi(j+1)}{\Phi(j+r)}
\bigl(\Phi (n) \bigr)^{r-1}.
\]
\item[2)] Suppose $r\in(0,1]$. Then
\[
\frac{h(n+1,r)}{h(n,r)}=\frac{\Phi(n+1)^r \Phi(n)^{1-r}}{\Phi
(n+r)}\leq1
\]
since, as $1/\Phi$ is completely monotone, $1/\Phi$ is log-convex and
therefore, $\log\Phi$ is concave.
Thus, $ n\longrightarrow h(n,r)$ is decreasing and, as it was seen in
the proof of Theorem~\ref{t21}, $h(n,r)\geq R(r)$. Therefore,
$ {h(r):=
\lim_{n\uparrow\infty}\downarrow h(n,r)}$ exists for $r\in(0,1]$ and
is $>0$.
\item[3)] One has, for $n\geq1$ and $r>0$,
\[
h(n,r+1)=\Phi(r) \biggl(\frac{\Phi(n)}{\Phi(n+1)} \biggr)^r h(n+1,r).
\]
Hence, taking into account the above step 2), we obtain by induction:
\[
\forall r>0 ,\quad h(r):=\li{n} {\infty}h(n,r) \mbox{ exists, and } h(r+1)=
\Phi(r) h(r).
\]
\item[4)] As $1/\Phi$ is log-convex, $h$ is log-convex too. Thus, $h$
satisfies the conditions i), ii) and iii) in Theorem \ref{t21}, and
consequently, $h=R$.\quad\qed
\end{enumerate}
\noqed\end{pf}
%
\begin{thmm}\label{t24} There is the asymptotic formula:
\[
\forall r>0 ,\quad I(r)=\Gamma(r)\li{n} {\infty} \prod
_{j=0}^{n-1}\frac
{\Phi(j+r)}{\Phi(j+1)} \bigl(\Phi(n)
\bigr)^{1-r}.
\]
\end{thmm}
\begin{pf}
This is a direct consequence of Theorem \ref{t23} and (\ref{e23}).
Another proof consists in using the same proof as that of Theorem \ref
{t23} (replacing $1/\Phi(s)$ by $\Phi(s)/s$) and (\ref{e24}).
\end{pf}
%
\subsection{Convolution equations satisfied by the laws of $\mathcal
{I}$ and $\mathcal{R}$}\label{ss25}

In this subsection, we assume, for simplicity, $a=0$. We also denote by
$\overline{\l}$ the density of the measure $\overline{\l}$:
$\overline
{\l}(x)=\l((x,+\infty))$ for $x\geq0$.

Property 1 in the following theorem is a particular case (Example B) of
Carmona et~al. (\cite{CPY}, Proposition 2.1). However, the proofs are different.
%
\begin{thmm}\label{t25}
1. The law of $\mathcal{I}$ admits a density which we denote by
$\PR
{\theta(v)}{v>0}$ and which satisfies:
%
\begin{equation}
\label{e25} \forall v>0,\quad\theta(v)=\Io{v} {\infty} {\theta(y) \overline{\l }
\biggl(\log \biggl(\frac{y}{v} \biggr) \biggr)} {y}.
\end{equation}
\begin{enumerate}[2.]
\item[2.] The law of $\mathcal{R}$ satisfies:
%
\begin{equation}
\label{e26} \forall v>0,\quad\mathbb{P}(\mathcal{R}>v)=\mathbb{E} \biggl[
\frac
{1}{\mathcal{R}} 1_{(\mathcal{R}>v)} \overline{\l} \biggl(\log \biggl(
\frac{\mathcal{R}}{v} \biggr) \biggr) \biggr].
\end{equation}
\item[3.] If the potential measure $\rho$ admits a density (still denoted
by $\rho$), then the law of $\mathcal{R}$ admits a density which we
denote by $\PR{\zeta(v)}{v>0}$ and which satisfies:
%
\begin{equation}
\label{e27} \forall v>0,\quad\zeta(v)=\Io{v} {\infty} {\zeta(y) \rho \biggl(\log
\biggl(\frac{y}{v} \biggr) \biggr)} {y}.
\end{equation}
\end{enumerate}
\end{thmm}
\begin{pf}
\begin{enumerate}[1)]
\item[1)] Let $\theta$ be the law of $\mathcal{I}$. By Proposition \ref
{p21}, we have, for $r>0$,
\[
\frac{\Phi(r)}{r} I(r+1)=\I{} {\int\mathrm{e}^{-xr} y^r
\overline {\l }(x) \,\mathrm{d}x} {\theta} {y}.
\]
The change of variable (from $x$ to $v$): $x=\log (\frac
{y}{v}
)$ yields:
\[
\frac{\Phi(r)}{r} I(r+1)=\Io{0} {\infty} {v^{r-1} \biggl[\int
_v^{\infty
}\overline{\l} \biggl(\log \biggl(
\frac{y}{v} \biggr) \biggr) \theta (\mathrm{d}y) \biggr]} {v}.
\]
Now,
\[
I(r)=\I{} {v^{r-1}} {\theta} {v}.
\]
Then, from the functional equation (\ref{e21}), the measures:
\[
\theta(\mathrm{d}v)\quad\mbox{and}\quad \biggl[\int_v^{\infty
}
\overline {\l} \biggl(\log \biggl(\frac{y}{v} \biggr) \biggr) \theta (
\mathrm{d}y) \biggr] \,\mathrm{d}v
\]
admit the same Mellin transform, hence they are equal.
\item[2)] We set $\eta(y)=\mathbb{P}(\mathcal{R}>y)$. We have, for $r>0$,
\[
\frac{\Phi(r)}{r} R(r)=-\int\int\mathrm{e}^{-xr} y^{r-1}
\overline {\l }(x) \,\mathrm{d}x \,\mathrm{d}\eta(y).
\]
The change of variable (from $x$ to $v$): $x=\log (\frac
{y}{v}
)$ yields:
\[
\frac{\Phi(r)}{r} R(r)=-\Io{0} {\infty} {v^{r-1} \biggl[\int
_v^{\infty
}\frac{1}{y} \overline{\l} \biggl(
\log \biggl(\frac{y}{v} \biggr) \biggr) \,\mathrm{d}\eta (y) \biggr]} {v}.
\]
Now,
\[
\frac{1}{r} R(r+1)=-\frac{1}{r} \int v^r \,\mathrm{d}
\eta(v)=\Io {0} {\infty} { v^{r-1} \eta(v)} {v}.
\]
Then (\ref{e26}) follows from the functional equation (\ref{e22}) by
injectivity of the Mellin transform.
\item[3)] The proof of point 3 is quite similar to that of point
1.\quad\qed
\end{enumerate}
\noqed\end{pf}
%
\subsection{The symmetry case}\label{ss26}
%
\begin{Def}\label{d22}
For any Bernstein function $\Phi\not\equiv0$, we define $\Phi^{\ast
}$ by:
\[
\forall s>0,\ \Phi^{\ast}(s)=\frac{s}{\Phi(s)}\quad\mbox{and}\quad
\Phi^{\ast}(0)=\li{s} {0}\Phi^{\ast}(s).
\]
We denote by $\Sigma$ the set of Bernstein functions $\Phi\not\equiv0$
such that $\Phi^{\ast}$ is a Bernstein function. (In particular (see
Definition \ref{d21}), if $\Phi\in\Sigma$, then $\Phi^{\ast
}(0)=0$.)

In the terminology of Schilling et al. \cite{SSV}, $\Sigma$ is the set
of \emph{special Bernstein functions} $\Phi$ such that
\[
\Phi(0)=0\quad\mbox{and}\quad\li{s} {0}\frac{\Phi
(s)}{s}=+\infty.
\]
We call \emph{symmetry case} the situation where $\Phi\in\Sigma$. In
this case:
\[
\forall s\geq0,\quad\Phi(s) \Phi^{\ast}(s)=s.
\]
Obviously, if $\Phi\in\Sigma$, then $\Phi^{\ast}\in\Sigma$ and
$\Phi^{\ast\ast}=\Phi$. We then say that $(\Phi,\Phi^{\ast})$ is a
\emph{pair
of conjugate Bernstein functions}.
\end{Def}
Our interest in the symmetry case stems from the next proposition.
%
\begin{pro}\label{p23}
Suppose that $(\Phi,\Phi^{\ast})$ is a pair of conjugate Bernstein
functions. We denote by $\mathcal{I}_{\Phi}$ and $\mathcal{R}_{\Phi}$,
respectively the perpetuity and the remainder variables related to the
L\'evy process whose Laplace--Bernstein exponent is $\Phi$, and likewise
for $\mathcal{I}_{\Phi^{\ast}}$ and $\mathcal{R}_{\Phi^{\ast}}$. Then:
\[
\mathcal{I}_{\Phi}\st\mathcal{R}_{\Phi^{\ast}}\quad\mbox{and}\quad
\mathcal{R}_{\Phi}\st\mathcal{I}_{\Phi^{\ast}}.
\]
\end{pro}
\begin{pf}
By Theorem \ref{t21} and Theorem \ref{t22}, $\mathcal{I}_{\Phi}$ and
$\mathcal{R}_{\Phi^{\ast}}$ have the same Mellin transform. Likewise,
$\mathcal{R}_{\Phi}$ and $\mathcal{I}_{\Phi^{\ast}}$ have the same
Mellin transform. The desired result follows from the injectivity of
the Mellin transform.
\end{pf}
We now state a useful characterization of the symmetry case (see also
Schilling et al. (\cite{SSV}, Theorem 10.3)).
%
\begin{thmm}\label{t26}
Let $\Phi\not\equiv0$ be a Bernstein function. Following Subsection
\ref{ss21}, we still denote by $\rho$ the related potential measure.
Then, $\Phi\in\Sigma$ if and only if there exist $b\geq0$ and
$h:(0,+\infty)\longrightarrow\R_+$ such that:
\begin{eqnarray}\label{e28}
&\displaystyle h \mbox{ is decreasing, } \li{x} {\infty}h(x)=0\quad\mbox{and}&
\nonumber
\\[-8pt]
\\[-8pt]
\nonumber
&\displaystyle \rho(\mathrm{d}x)=b \e(\mathrm{d}x)+h(x) \,\mathrm{d}x.&
\end{eqnarray}
Moreover, if these properties are satisfied, then:
%
\begin{equation}
\label{e29} \Phi^{\ast}(s)=b s-\int_0^{\infty}
\bigl(1-\mathrm{e}^{-sx}\bigr) \,\mathrm {d}h(x),\ s\geq0.
\end{equation}
\end{thmm}
\begin{pf}
\begin{enumerate}[1)]
\item[1)] Suppose $\Phi\in\Sigma$. Then
%
\begin{equation}
\label{e210} \Phi^{\ast}(s)=a^{\ast} s+\I{} {\bigl(1-
\mathrm{e}^{-sx}\bigr)} {\l^{\ast}} {x}.
\end{equation}
Since $ {\frac{1}{\Phi(s)}=\frac{\Phi^{\ast}(s)}{s}}$,
Proposition \ref{p21} yields:
\[
\rho(\mathrm{d}x)=a^{\ast} \e(\mathrm{d}x)+\l^{\ast}\bigl((x,+
\infty)\bigr) \,\mathrm{d}x.
\]
This shows the ``only if'' part, setting:
%
\begin{equation}
\label{e211} b=a^{\ast}\quad\mbox{ and }\quad h(x)=\l^{\ast}
\bigl((x,+\infty)\bigr).
\end{equation}
Moreover, (\ref{e210}) and (\ref{e211}) entail (\ref{e29}).
\item[2)] Conversely, suppose that (\ref{e28}) holds. Since $\rho$ is a
Radon measure on $\R_+$, then
\[
\Io{0} {1} {h(x)} {x}=\int_0^{\infty}(y\wedge1)
\bigl(-\mathrm {d}h(y)\bigr)<\infty.
\]
Therefore, the Stieltjes measure $-\mathrm{d}h$ is a L\'evy measure and
the function:
\[
\Psi(s)=b s-\int_0^{\infty}\bigl(1-
\mathrm{e}^{-sx}\bigr) \,\mathrm{d}h(x),\ s\geq0
\]
is a Bernstein function. Using again Proposition \ref{p21} and (\ref
{e28}), we obtain:
\[
\forall s>0,\quad\frac{\Psi(s)}{s}=\frac{1}{\Phi(s)}.
\]
Hence, $\Phi\in\Sigma$ and $\Phi^{\ast}=\Psi$.\quad\qed
\end{enumerate}
\noqed\end{pf}
Now, we exhibit a large convex cone $\Lambda$ contained in $\Sigma$
(see also Schilling et al. (\cite{SSV}, Theorem~10.11)).
%
\begin{thmm}\label{t27}
We denote by $\Lambda$ the convex cone consisting of functions $\Psi$
of the form:
\[
\Psi(s)=b s-\int_0^{\infty}\bigl(1-
\mathrm{e}^{-sx}\bigr) \,\mathrm{d}h(x),\ s\geq0,
\]
for some $b\geq0$ and some decreasing log-convex function $h$ on
$(0,+\infty)$ satisfying:
\[
\Io{0} {1} {h(x)} {x}<\infty,\ \Io{0} {\infty} {h(x)} {x}=+\infty\quad \mbox
{and}\quad\li{x} {\infty}h(x)=0.
\]
Then, $\Lambda\subset\Sigma$.
\end{thmm}
\begin{pf}
Let $b\geq0$ and let $h$ be a decreasing log-convex function on
$(0,+\infty)$ satisfying
$\Io{0}{1}{h(x)}{x}<\infty$ and $\Io{0}{\infty}{h(x)}{x}=+\infty$.
We set:
\[
\rho(\mathrm{d}x)=b \e(\mathrm{d}x)+h(x) \,\mathrm{d}x.
\]
By Hirsch (\cite{H}, Th\'eor\`eme 2) (see also It\^{o} \cite{I}), $\rho$
is the potential measure associated with a Bernstein function $\Phi$.
By Theorem \ref{t26}, if moreover $\li{x}{\infty}h(x)=0$, then $\Phi
\in
\Sigma$ and
\[
\Phi^{\ast}(s)=b s-\int_0^{\infty}\bigl(1-
\mathrm{e}^{-sx}\bigr) \,\mathrm {d}h(x),\ s\geq0.
\]
But, as we have noted before, $\Psi:=\Phi^{\ast}\in\Sigma$.
\end{pf}
%
\begin{Def}\label{d23}
Following Schilling et al. \cite{SSV}, we call \emph{complete Bernstein
function} any function $\Phi$ of the form:
%
\begin{equation}
\label{e212} \Phi(s)=b+a s+\int_0^{\infty}\bigl(1-
\mathrm{e}^{-sx}\bigr) m(x) \,\mathrm{d}x,\ s\geq0,
\end{equation}
for some $a\geq0$, $b\geq0$, and some completely monotone function
$m$ on $(0,+\infty)$ satisfying:
\[
\Io{0} {\infty} {(x\wedge1) m(x)} {x}<\infty.
\]
We denote by $\widehat{\mathcal{S}}$ the convex cone of complete
Bernstein functions.

We also denote by $\mathcal{S}$ the convex cone consisting of functions
$\Phi\in\widehat{\mathcal{S}}$ such that:
\[
\Phi(0)=0 \quad\mbox{and}\quad\Phi^{\ast}(0)= \li{s} {0}
\frac
{s}{\Phi(s)}=0,
\]
or, with the notation of (\ref{e212}),
%
\begin{equation}
\label{ne21} b=0\quad\mbox{and}\quad\Io{0} {\infty} {x m(x)} {x}=\infty.
\end{equation}
Using Bernstein's theorem, we obtain that $\Phi\in\mathcal{S}$ if and
only if there exist $a\geq0$ and a measure~$\mu$ on $(0,+\infty)$ so that:
\begingroup
\abovedisplayskip=6.5pt
\belowdisplayskip=6.5pt
\begin{eqnarray*}
&\displaystyle \I{(0,+\infty)} {\frac{1}{1+t}} {\mu} {t}<\infty,\ \I{(0,+\infty )} {
\frac
{1}{t}} {\mu} {t}=+\infty\quad \mbox{and}&
\\
&\displaystyle\Phi(s)=a s+\I{(0,+\infty)} {\frac{s}{s+t}} {\mu} {t},\ s\geq0.&
\end{eqnarray*}
In other words, $\Phi\in\mathcal{S}$ if and only if $ {\frac
{\Phi(s)}{s}}$ is a \emph{Stieltjes transform} $S$ such that
\[
\li{s} {0}S(s)=+\infty\quad\mbox{and}\quad\li{s} {0}s S(s)=0.
\]

\end{Def}

\begin{pro}\label{p24}
There is the double inclusion: $\mathcal{S}\subset\Lambda\subset
\Sigma
$. Moreover, $\Phi$ belongs to $\mathcal{S}$ if and only if $\Phi^{\ast
}$ belongs to $\mathcal{S}$.
\end{pro}
\begin{pf}
\begin{enumerate}[1)]
\item[1)] Suppose $\Phi\in\mathcal{S}$ and $\Phi$ given by (\ref
{e212}) and
(\ref{ne21}). Set:
\[
\forall x>0,\quad h(x)=\Io{x} {\infty} {m(y)} {y}.
\]
Then,
\begin{eqnarray*}
&\displaystyle \Io{0} {1} {h(x)} {x}=\Io{0} {\infty} {(x\wedge1) m(x)} {x}<\infty\quad\mbox{and}&
\\
&\displaystyle \Io{0} {\infty} {h(x)} {x}=\Io{0} {\infty} {x m(x)} {x}=+\infty.&
\end{eqnarray*}

Moreover, $\li{x}{\infty}h(x)=0$, $h$ is decreasing, and, since $h$ is
completely monotone, then $h$ is log-convex. Therefore, $\Phi\in
\Lambda$.
\item[2)] By Schilling et al. (\cite{SSV}, Proposition 7.1), if $\Phi\not
\equiv0$ is a complete Bernstein function, then so is $\Phi^{\ast}$.
This gives another proof of $\mathcal{S}\subset\Sigma$ and shows that,
if $\Phi\in\mathcal{S}$, then $\Phi^{\ast}\in\mathcal{S}$ too.\hfill\qed
\end{enumerate}
\noqed\end{pf}\endgroup

The cone $\widehat{\mathcal{S}}$ of complete Bernstein functions has a
deep probabilistic interpretation. Indeed, using Krein's theory of
strings, an open problem in It\^{o}--Mc Kean concerning the precise
class of subordinators which are inverse local times of a regular
diffusion on $[0,+\infty)$, instantaneously reflecting at 0, has been
solved simultaneously and independently in 1981 by Kotani--Watanabe~\cite
{KW} and Knight \cite{K}. See also Bertoin \cite{Be}, K\"{u}chler~\cite{Ku},
K\"{u}chler--Salminen \cite{KS}. In fact, this class is precisely
the class of subordinators whose Laplace--Bernstein exponent belongs to
$\widehat{\mathcal{S}}$. Interestingly, in Bertoin's Saint-Flour course
\cite{Be2}, Chapter 9, there is the description of $\Phi^{\ast}$, where
$\Phi$ is the Bernstein function associated to:
\[
A_{\tau_l}:=\Io{0} {\tau_l} {f(B_t)} {t} ,\quad l
\geq0,
\]
with $f\geq0$ with support in $(0,+\infty)$, and $\PR{\tau_l}{l\geq
0}$ the inverse local time at 0 of the Brownian motion $\PR{B_t}{t\geq
0}$. Precisely, Corollary 9.7, page 78, asserts that in that case:
$\Phi^{\ast}\in\widehat{\mathcal{S}}$. More generally, the discussion of
Krein's theory of strings and its relationship with \emph{inverse local
times of generalized diffusions} is expounded in an exhaustive manner
in Chapter 14 of Schilling et al. \cite{SSV}. The Krein table (that is:
identifying the generalized diffusion whose Laplace--Bernstein exponent
of the inverse local time is a given element of $\widehat{\mathcal
{S}}$) established by Donati-Martin and Yor (\cite{DY} and \cite{DY2})
is reproduced in Schilling et al. \cite{SSV} on pp. 201--202. Note that
it has relatively few entries, and would certainly deserve to be completed.
\section{Multiplicative infinite divisibility}\label{s3}
\subsection{Multiplicative infinite divisibility of $\mathcal
{R}$}\label{ss31}
%
\begin{Def}\label{d31}
A positive random variable $X$ is said to be \emph{multiplicatively
infinitely divisible} (m.i.d.) if, for any $n\in\N$ with $n\geq2$,
there exist independent identically distributed positive random
variables: $X_1,\ldots,X_n$, such that:
\[
X\st X_1\cdots X_n.
\]
\end{Def}
Obviously, a strictly positive random variable $X$ is m.i.d. if and
only if $\log X$ is infinitely divisible.

In the sequel, the notation of the previous section is still in force.
The following proposition was also proven by C. Berg (\cite{B3}, Theorem 1.8).
%
\begin{pro}\label{p31}
The remainder $\mathcal{R}$ is m.i.d.
\end{pro}
\begin{pf}
For $n\geq2$, let $\Phi_n:=(\Phi)^{1/n}$. As it is well known, $\Phi_n$ also is
a Bernstein function (associated with the subordinate L\'
evy process: $\xi^{(1/n)}:=\PR{\xi_{\tau^{(1/n)}_l}}{l\geq0}$ where
$\tau^{(1/n)}$ denotes a $(1/n)$-stable subordinator independent of
$\xi
$). Let $\mathcal{R}_n$ be the remainder related to $\xi^{(1/n)}$.
Then, we deduce easily from Theorem \ref{t21} and from the injectivity
of the Mellin transform:
\[
\mathcal{R}\st\mathcal{R}^{(1)}_n\cdots
\mathcal{R}^{(n)}_n
\]
where $\mathcal{R}^{(1)}_n,\ldots,\mathcal{R}^{(n)}_n$ are $n$
independent copies of
$\mathcal{R}_n$. Thus, by Definition \ref{d31}, $\mathcal{R}$ is m.i.d.
\end{pf}
%
\subsection{Integral representations}\label{ss32}
We shall deduce from Proposition \ref{p31} a representation of the
Mellin transform $R$ of $\mathcal{R}$. We recall (cf. Subsection \ref
{ss21}) that $\kappa$ denotes the measure whose Laplace transform is
$\Phi'/\Phi$.
%
\begin{thmm}\label{t31}
We have, for $r>0$,
\[
R(r)=\Phi(1)^{r-1} \exp \biggl[\I{(0,+\infty)} {\frac{\mathrm
{e}^{-(r-1)x}-1-(r-1)(\mathrm{e}^{-x}-1)}{x (\mathrm{e}^x-1)}} {
\kappa } {x} \biggr].
\]
\end{thmm}
\begin{pf}
We shall give two proofs. We mention that C. Berg stated a more general
result (see Berg (\cite{B2}, Theorem 2.2)), with a different proof.

\textit{First proof}. By Proposition \ref{p31} and the L\'evy--Khintchine formula, there exist $a,\sigma\in\R$ and a measure
$\Theta$ on $\R\setminus\{0\}$ satisfying
\[
\I{} {\frac{x^2}{1+x^2}} {\Theta} {x}<\infty
\]
such that
\[
\forall u\in\R,\quad\mathbb{E} \bigl[\mathcal{R}^{\mathrm
{i}u} \bigr]=\exp
\bigl(-\psi(u)\bigr)
\]
with
\[
\psi(u)=\mathrm{i} a u+\frac{1}{2} \sigma^2 u^2+
\I{} {\bigl(1-\mathrm {e}^{\mathrm{i}ux}+\mathrm{i} u x 1_{\{|x|<1\}}\bigr)} {
\Theta} {x}.
\]
Since the Mellin transform $R(r)$ is defined for $r>0$, then $\psi$
continuously extends to $\mathrm{i} (-\infty, 1)$. Hence,
\[
\forall s>0,\ \I{(1,+\infty)} {\mathrm{e}^{sx}} {\Theta} {x}<\infty
\quad \mbox{and}\quad\forall0<s<1,\ \I{(-\infty,-1)} {\mathrm {e}^{-sx}} {
\Theta} {x}<\infty
\]
and, for $r>0$,
%
\begin{equation}
\label{e31} \log\bigl(R(r)\bigr)=-a (r-1)+\frac{1}{2} \sigma^2
(r-1)^2+\I{} {\bigl(\mathrm {e}^{(r-1)x}-1-(r-1) x
1_{\{|x|<1\}}\bigr)} {\Theta} {x}.
\end{equation}
We deduce then from (\ref{e31}) and (\ref{e22}):
\[
\log\Phi(r)=-a+\frac{1}{2} \sigma^2 (2r-1)+\I{} {\bigl(
\mathrm {e}^{rx}-\mathrm{e}^{(r-1)x}-x 1_{\{|x|<1\}}\bigr)} {
\Theta} {x},
\]
and, by differentiation,
\[
\frac{\Phi'(r)}{\Phi(r)}=\sigma^2+\I{} {x \bigl(\mathrm{e}^{rx}-
\mathrm {e}^{(r-1)x}\bigr)} {\Theta} {x}.
\]
Therefore, $\sigma=0$ and the measure $x (1-\mathrm{e}^{-x}) \Theta
(\,\mathrm{d}x)$ is the image of $\kappa$ by $x\longrightarrow-x$. (In
particular, $\Theta$ is carried by $\R_-$.)

Then, (\ref{e31}) becomes:
%
\begin{equation}
\label{e32} \log\bigl(R(r)\bigr)=-a (r-1)+\I{} {\frac{(\mathrm{e}^{-(r-1)x}-1+(r-1) x 1_{\{
0<x<1\}})}{x (\mathrm{e}^x-1)}} {\kappa} {x}
\end{equation}
and, in particular, for $r=2$,
%
\begin{equation}
\label{e33} \log\bigl(\Phi(1)\bigr)=-a+\I{} {\frac{(\mathrm{e}^{-x}-1+x 1_{\{0<x<1\}})}{x
(\mathrm{e}^x-1)}} {\kappa} {x}.
\end{equation}
The desired result follows directly from (\ref{e32}) and (\ref{e33}).

\textit{Second proof}. We set, for $r>0$,
\[
f(r)=\Phi(1)^{r-1} \exp \biggl[\I{(0,+\infty)} {\frac{\mathrm
{e}^{-(r-1)x}-1-(r-1)(\mathrm{e}^{-x}-1)}{x (\mathrm{e}^x-1)}} {
\kappa } {x} \biggr].
\]
Obviously, $f$ is a log-convex function such that $f(1)=1$. Therefore,
by Theorem \ref{t21}, we only have to prove:
\[
\forall r>0, \quad f(r+1)=\Phi(r) f(r).
\]
Now, a simple computation yields:
\begin{eqnarray*}
\log f(r+1)-\log f(r)&=&\log\Phi(1)+\I{} {\frac{\mathrm
{e}^{-x}-\mathrm
{e}^{-rx}}{x}} {\kappa} {x}
\\
&=&\log\Phi(1)+\Io{1} {r} { \biggl(\I{} {\mathrm{e}^{-sx}} {\kappa } {x}
\biggr)} {s}
\\
&=&\log\Phi(1)+\Io{1} {r} {\frac{\Phi'(s)}{\Phi(s)}} {s}=\log\Phi(r).
\end{eqnarray*}
\upqed\end{pf}

In the trivial example, Theorem \ref{t31} yields a classical
representation of the Gamma function (see, for instance, Andrews et al.
(\cite{AAR}, Theorem 1.6.2)):
%
\begin{equation}
\label{e34} \Gamma(r)=\exp \biggl[\Io{0} {+\infty} {\frac{\mathrm
{e}^{-(r-1)x}-1-(r-1)(\mathrm{e}^{-x}-1)}{x (\mathrm
{e}^x-1)}} {x}
\biggr].
\end{equation}
We now state a result which was seen in the first proof of Theorem \ref{t31}.
%
\begin{pro}\label{p32}
The L\'evy measure of the infinite divisible variable $\log\mathcal{R}$
is the image by the map $x\longrightarrow-x$ of the measure on
$(0,+\infty)$:
\[
x^{-1} \bigl(\mathrm{e}^{x}-1\bigr)^{-1} \kappa(
\mathrm{d}x).
\]
\end{pro}
As a direct consequence of Theorem \ref{t31} and (\ref{e23}), we obtain
the following theorem.
%
\begin{thmm}\label{t32}
For $r>0$,
\[
I(r)=\Gamma(r) \Phi(1)^{-r+1} \exp \biggl[-\I{(0,+\infty)} {
\frac
{\mathrm
{e}^{-(r-1)x}-1-(r-1)(\mathrm{e}^{-x}-1)}{x (\mathrm{e}^x-1)}} {\kappa } {x} \biggr].
\]
\end{thmm}
Note that, by (\ref{e34}), we also have, for $r>0$,
%
\begin{equation}
\label{e35} I(r)=\Phi(1)^{-r+1} \exp \biggl[\int_{(0,+\infty)}
\frac{\mathrm
{e}^{-(r-1)x}-1-(r-1)(\mathrm{e}^{-x}-1)}{x (\mathrm{e}^x-1)} \bigl(\mathrm {d}x-\kappa(\mathrm{d}x)\bigr) \biggr].
\end{equation}
%
\subsection{Multiplicative infinite divisibility of $\mathcal
{I}$}\label{ss33}
In the next theorem, we state a characterization of the multiplicative
infinite divisibility of $\mathcal{I}$ (see also Berg (\cite{B3}, Theorem 1.9)).
%
\begin{thmm}\label{t33}
The variable $\mathcal{I}$ is m.i.d. if and only if $\kappa(\mathrm
{d}x)\leq\mathrm{d}x$. Besides, if $\mathcal{I}$ is m.i.d., then the
L\'
evy measure of the infinitely divisible variable $\log\mathcal{I}$ is:
\[
1_{(-\infty,0)}(x) x^{-1} \bigl(1-\mathrm{e}^{-x}
\bigr)^{-1} \bigl(1-k(-x)\bigr) \,\mathrm{d}x
\]
where $k$ denotes the density: $ {\frac{\kappa(\mathrm
{d}x)}{\mathrm{d}x}}$.
\end{thmm}

\begin{pf}
We deduce from formula (\ref{e35}) that:
\[
\forall u\in\R,\quad\mathbb{E} \bigl[\mathcal{I}^{\mathrm
{i}u} \bigr]=\exp
\bigl(-\eta(u)\bigr)
\]
with
\[
\eta(u)=\mathrm{i} u \log\Phi(1)+\int_{(0,+\infty)}\frac
{1-\mathrm
{e}^{-\mathrm{i}ux}+\mathrm{i} u (\mathrm{e}^{-x}-1)}{x (\mathrm
{e}^x-1)}
\bigl(\mathrm{d}x-\kappa(\mathrm{d}x)\bigr).
\]
Then the result follows from the L\'evy--Khintchine formula.
\end{pf}
We now give a sufficient condition.
%
\begin{pro}\label{p33}
If $\Phi\in\Sigma$ (see Definition \ref{d22}), then the variable
$\mathcal{I}$ is m.i.d.
\end{pro}
\begin{pf}
If $\Phi\in\Sigma$, then, with obvious notation,
\[
\kappa(\mathrm{d}x)+\kappa^{\ast}(\mathrm{d}x)=1_{\R_+}(x)
\,\mathrm{d}x.
\]
In particular, $\kappa(\mathrm{d}x)\leq\mathrm{d}x$, and Theorem
\ref{t33} applies.

Another proof consists in using jointly Proposition \ref{p23} and
Proposition \ref{p31}.
\end{pf}
Another sufficient condition is the following proposition.
%
\begin{pro}\label{p34}
Suppose that, for every $\a\in(0,1)$, the function $s^{1-\a} \Phi^{\a
}(s)$ is a Bernstein function, then the variable $\mathcal{I}$ is m.i.d.
\end{pro}
\begin{pf}
For $n\in\N$, $n\geq2$, we denote by $\Phi_n$ the Bernstein function
$s^{1-\a} \Phi^{\a}(s)$ with $\a=1/n$. Then
\[
\frac{s}{\Phi(s)}= \biggl(\frac{s}{\Phi_n(s)} \biggr)^n.
\]
We denote by $\mathcal{I}_n$ the perpetuity related to the Bernstein
function $\Phi_n$. Then, we deduce easily from Theorem \ref{t22} and
from the injectivity of the Mellin transform:
\[
\mathcal{I}\st\mathcal{I}^{(1)}_n\cdots
\mathcal{I}^{(n)}_n
\]
where $\mathcal{I}^{(1)}_n\cdots\mathcal{I}^{(n)}_n$ are $n$
independent copies of
$\mathcal{I}_n$. Thus, by Definition \ref{d31}, $\mathcal{I}$ is m.i.d.
\end{pf}
%
\begin{cor}\label{c31}
Suppose that $\Phi$ is a complete Bernstein function (see Definition
\ref{d23}), null at 0.
Then, the variable $\mathcal{I}$ is m.i.d.
\end{cor}
\begin{pf}
By Schilling et al. (\cite{SSV}, Proposition 7.10.) (see also Berg \cite
{B}), a complete Bernstein function $\Phi$ satisfies the condition of
Proposition \ref{p34}, which entails the desired result.

Note that if moreover $\Phi^{\ast}(0)=0$, then $\Phi\in\mathcal{S}$,
and the multiplicative infinite divisibility of $\mathcal{I}$ also
follows from Proposition \ref{p24} and Proposition \ref{p33}.
\end{pf}

We end this section by a straightforward consequence of Proposition
\ref
{p32} and Theorem~\ref{t33}.
%
\begin{pro}\label{p35}
The variable $\log\mathcal{R}$ is self-decomposable if and only if
there exists a positive decreasing function $\jmath$ on $ (0,+\infty)$
such that:
\[
\kappa(\mathrm{d}x)=1_{(0,+\infty)}(x) \bigl(\mathrm{e}^x-1\bigr)
\jmath(x) \,\mathrm{d}x.
\]
The variable $\log\mathcal{I}$ is self-decomposable if and only if
there exists a positive decreasing function~$\ell$ on $ (0,+\infty)$
such that:
\[
\kappa(\mathrm{d}x)=1_{(0,+\infty)}(x) \bigl[1-\bigl(\mathrm{e}^x-1
\bigr) \ell(x)\bigr] \,\mathrm{d}x.
\]

\end{pro}
%
\section{Examples}\label{s4}
\subsection{\texorpdfstring{$\a$-stable subordinator}{alpha-stable subordinator}}\label{ss41}
Let $\a\in(0,1)$ and let $\xi$ be an $\a$-stable subordinator.
Then, $\Phi(s)=s^{\a}$. Consequently, $\Phi\in\mathcal{S}$ and
$\Phi^{\ast}(s)=s^{1-\a}$. We have: $\kappa(\mathrm{d}x)=\a1_{\R_+}(x)
\,\mathrm{d}x$. Therefore,
$\log\mathcal{I}$ and $\log\mathcal{R}$ are self-decomposable.
Moreover, we obtain easily:
\[
I(r)=\bigl[\Gamma(r)\bigr]^{1-\a}\quad\mbox{and}\quad R(r)=\bigl[
\Gamma(r)\bigr]^{\a}.
\]

The diffusion whose inverse local time is an $\a$-stable subordinator
is a Bessel process of dimension $d=2(1-\a)$ (see, for instance,
Molchanov--Ostrovski \cite{MO}).
\subsection{Exponential compound Poisson process}\label{ss42}
Let $c>0$ and let $\xi$ be a compound Poisson process whose L\'evy
measure is
$c 1_{\R_+}(x) \mathrm{e}^{-cx} \,\mathrm{d}x$. Then, $\Phi(s)=s
(s+c)^{-1}$. Consequently, $\Phi$ is a complete Bernstein function
($m(x)=c \mathrm{e}^{-cx}$), but $\Phi\notin\mathcal{S}$ (since
$\Phi^{\ast}(s)=s+c$ and $\Phi^{\ast}(0)=c\not=0$). We have: $\kappa
(\mathrm
{d}x)=1_{\R_+}(x) (1-\mathrm{e}^{-cx}) \,\mathrm{d}x$. Thus $\mathcal
{I}$ is m.i.d. More precisely, we deduce from Proposition \ref{p35}
that the variables $\log\mathcal{I}$ and $\log\mathcal{R}$ are
self-decomposable. Moreover, we obtain easily (for example by Theorem
\ref{t25}):
\[
I(r)=\frac{\Gamma(c+r)}{\Gamma(c+1)}\quad\mbox{and}\quad R(r)=\frac
{\Gamma(c+1) \Gamma(r)}{\Gamma(c+r)}=c B(r,c).
\]
Therefore, $\mathcal{I}\st\gamma_{c+1}$ and $ \mathcal{R}\st\beta_{1,c}$ where $\gamma_u$ denotes a gamma variable of parameter $u$ and
$\beta_{u,v}$ denotes a beta variable of parameters $u,v$.
\subsection{Geometric compound Poisson process}\label{ss45}
Let $0\leq c<q<1$. Following Bertoin et al. \cite{BBY}, we consider a
compound Poisson process $\PR{\xi_l}{l\geq0}$ whose L\'evy measure is
the geometric probability:
\[
\l=\sum_{n=1}^{\infty}(c/q)^{n-1}
(1-c/q) \e_{-n\log q}
\]
where $\e_x$ denotes the Dirac measure at point $ x$. (In particular,
if $c=0$, then $\l=\e_{-\log q}$ and $\xi_l=-(\log q) N_l$ where
$\PR
{N_l}{l\geq0}$ denotes the standard Poisson process.)

We obtain easily:
\[
\Phi(s)=\frac{1-q^s}{1-c q^{s-1}}.
\]
Then, by Theorems \ref{t23} and \ref{t24}, we have:
\begin{eqnarray*}
&\displaystyle R(r)=\prod_{j=0}^{\infty}\frac{(1-q^{j+1}) (1-c
q^{j+r-1})}{(1-q^{j+r}) (1-c q^j)}&
\\
&\displaystyle\mbox{and}\quad I(r)=\Gamma(r) \prod_{j=0}^{\infty}
\frac
{(1-q^{j+r}) (1-c q^j)}{(1-q^{j+1}) (1-c q^{j+r-1})}.&
\end{eqnarray*}
These formulae are proven in Bertoin et al. \cite{BBY}, where they are
related to the so-called $q$-calculus (see, for instance, Gasper--Rahman
\cite{GR}).

It is easy to see that the measure $\kappa$ is given by:
\[
\kappa=-\log q \sum_{n=1}^{\infty}
\bigl(1-(c/q)^n\bigr) \e_{-n\log q}.
\]
Hence, by Theorem \ref{t33}, $\mathcal{I}$ is not m.i.d. However, it is
proven in Bertoin et al. \cite{BBY} that the variable~$\mathcal{I}_0$,
corresponding to $c=0$, i.e.:
\[
\mathcal{I}_0=\Io{0} {\infty} {q^{N_l}} {l},
\]
is self-decomposable. Moreover, one has:
\[
\mathcal{I}_0 \st c \mathcal{I}_0+\mathcal{I}
\]
where, in the RHS, the variables are assumed to be independent (see
Bertoin et al. (\cite{BBY}, Proposition 3.1)). This entails that the
perpetuity $\mathcal{I}$ is infinitely divisible.

\subsection{Gamma process}\label{ss46}
We assume here that $\PR{\xi_l}{l\geq0}$ is the Gamma process. Thus,
\[
\forall l\geq0,\quad\xi_l\st\gamma_l\quad\mbox{and}
\quad\Phi (s)=\log(1+s).
\]
In particular, $\Phi$ is a complete Bernstein function ($m(x)=\mathrm
{e}^{-x}/x$), but, since $\Phi^{\ast}(0)=1$, $\Phi\notin\mathcal{S}$.
Nevertheless, by Corollary \ref{c31}, $\mathcal{I}$ is m.i.d.

We now determine the measure $\kappa$. An easy Laplace transform
computation (see also Remark~\ref{r51} and formula (\ref{e57}) below) shows
that the density $k$ of $\kappa$ is given by:
%
\begin{equation}
\label{ne51} k(x)=\mathrm{e}^{-x} \Io{0} {\infty} {\frac{x^l}{\Gamma(l+1)}}
{l},
\end{equation}
which entails:
%
\begin{equation}
\label{ne52} k'(x)=\mathrm{e}^{-x} \Io{0} {1} {
\frac{x^{l-1}}{\Gamma(l)}} {l},
\end{equation}
and therefore:
\[
k(x)=\Io{0} {x} {\Io{0} {1} {\frac{\mathrm{e}^{-y} y^{l-1}}{\Gamma
(l)}} {l}} {y}=\mathbb{P}(
\gamma_U\leq x)
\]
where $U$ is uniform on $(0,1)$ and independent from $\PR{\gamma_l}{l\geq0}$.
Consequently, $k$ is an increasing function, and $\li{x}{\infty}
k(x)=1$. In particular, the function: $\ell(x)=(\mathrm{e}^x-1)^{-1}
(1-k(x))$ is decreasing and hence, by Proposition \ref{p35}, $\log
\mathcal{I}$ is self-decomposable. Let: $\jmath(x)=(\mathrm
{e}^x-1)^{-1} k(x)$. We deduce easily from (\ref{ne51}) and (\ref
{ne52}) that:
\[
\jmath'(x)\leq\bigl(\mathrm{e}^x-1
\bigr)^{-2} \bigl(1-\mathrm{e}^{-x}-x\bigr) \Io {0} {1} {
\frac{x^{l-1}}{\Gamma(l)}} {l}\leq0
\]
and hence, by Proposition \ref{p35}, $\log\mathcal{R}$ is self-decomposable.

The diffusion whose inverse local time is a Gamma process was
determined by Donati-Martin and Yor \cite{DY}.
\subsection{Some examples from Bertoin--Yor}\label{ss43}
The next examples appear in Bertoin--Yor \cite{BY}.
\subsubsection{}\label{sss432}
Let $\a\in(0,1)$, $c>1$ and
\[
\Phi(s)=\frac{\a s \Gamma(\a( s-1+c))}{\Gamma(\a(s+c))}.
\]
By Bertoin--Yor \cite{BY}, $\Phi$ is a Bernstein function such that
\[
a=0,\quad\l(\mathrm{d}x)=m(x) \,\mathrm{d}x\quad\mbox{and}\quad
m(x)=-h'(x)
\]
with
\[
h(x)=\frac{1}{\Gamma(\a)} \frac{\mathrm{e}^{-(c-1)x}}{(1-\mathrm
{e}^{-x/\a})^{1-\a}}.
\]
Then, $h$ is a completely monotone function, hence $m$ is a completely
monotone function and $\Phi$ is a complete Bernstein function.
Consequently, by Corollary \ref{c31}, $\mathcal{I}$ is m.i.d.

We now determine the measure $\kappa$. We first remark that formula
(\ref{e34}) yields, by differentiation, the following classical formula:
%
\begin{equation}
\label{e41} \frac{\Gamma'}{\Gamma}(r)=\Io{0} {\infty} { \biggl(\frac{\mathrm
{e}^{-x}}{x}-
\frac{\mathrm{e}^{-rx}}{1-\mathrm{e}^{-x}} \biggr)} {x}.
\end{equation}
We deduce therefrom, by a simple computation,
\[
\kappa(\mathrm{d}x)=\frac{1-\mathrm{e}^{-\frac{1}{\a}x}+\mathrm
{e}^{-cx}-\mathrm{e}^{-(c-1)x}}{1-\mathrm{e}^{-\frac{1}{\a}x}} \,\mathrm{d}x.
\]
Then, with the notation of Proposition \ref{p35}, one has:
\[
\ell(x)=\frac{\mathrm{e}^{-cx}}{1-\mathrm{e}^{-\frac{1}{\a}x}}.
\]
Clearly, $\ell$ is a decreasing function, therefore, by Proposition
\ref
{p35}, $\log\mathcal{I}$ is self-decomposable. Moreover, we obtain
easily by Theorem \ref{t22} and (\ref{e23}):
\[
I(r)=\a^{1-r} \frac{\Gamma(\a(r-1+c))}{\Gamma(\a c)}\quad\mbox {and}\quad R(r)=
\a^{r-1} \frac{\Gamma(r) \Gamma(\a c)}{\Gamma(\a(r-1+c))}.
\]
Therefore, $\mathcal{I}\st\a^{-1} \gamma_{\a c}^{\a}$ and the law
of $
\mathcal{R}$ may also be made explicit (see Bertoin--Yor~\cite{BY}).

It may be noted that, in the case $\a=1/2$, the subordinator whose
Laplace--Bernstein exponent is the above complete Bernstein function
$\Phi$, appears in Comtet et al. (\cite{CTT}, Example 5.2).
\subsubsection{}\label{sss433}
Let $\a\in(0,1)$, $1< b\leq c$ and
\[
\Phi(s)=\frac{s \Gamma(\a( s+c))}{(b+s-1) \Gamma(\a(s-1+c))},\ \Phi^{\ast}(s)=\frac{(b+s-1) \Gamma(\a( s-1+c))}{\Gamma(\a(s+c))}.
\]
We have: $\Phi(0)=0$ and $\Phi^{\ast}(0)=(b-1) \Gamma(\a
(c-1))/\Gamma
(\a c)\not=0$. One can prove that:
\[
\Phi^{\ast}(s)=\Phi^{\ast}(0)-\Io{0} {\infty} {\bigl(1-
\mathrm{e}^{-sx}\bigr) \mathrm{e}^{-(b-1)x} h'(x)} {x}
\]
with
\[
h(x)=\frac{1}{\Gamma(1+\a)} \frac{\mathrm{e}^{-(c-b)x}}{(1-\mathrm
{e}^{-x/\a})^{1-\a}}.
\]
Then, $h$ is a completely monotone function, hence $\Phi^{\ast}$ is a
complete Bernstein function. By Schilling et al. (\cite{SSV}, Proposition 7.1),
$\Phi$ also is a complete Bernstein function. Consequently,
by Corollary \ref{c31}, $\mathcal{I}$ is m.i.d. It also may be seen that:
\[
\kappa(\mathrm{d}x)= \biggl[1- \biggl(\mathrm{e}^{-(b-1)x}-\mathrm
{e}^{-(c-1)x} \frac{1-\mathrm{e}^{-x}}{1-\mathrm{e}^{-\frac{1}{\a
}x}} \biggr) \biggr] \,\mathrm{d}x.
\]
Moreover, we obtain by Theorem \ref{t21} and (\ref{e23}):
\[
R(r)=\frac{b-1}{\Gamma(\a c)} B(r,b-1) \Gamma\bigl(\a(r-1+c)\bigr)\quad\mbox {and}
\quad I(r)=\frac{\Gamma(\a c) \Gamma(r-1+b)}{\Gamma(b) \Gamma(\a(r-1+c))}.
\]
Therefore, $\mathcal{R}\st\beta_{1,b-1} \gamma_{\a c}^{\a}$ with the
variables $\beta_{1,b-1}$ and $\gamma_{\a c}^{\a}$ independent, and the
law of $ \mathcal{I}$ may also be determined (see Bertoin--Yor \cite{BY}).
\subsection{Inverse local time of a radial Ornstein--Uhlenbeck
process}\label{ss44}
In this subsection, we assume that $\PR{\xi_l}{l\geq0}$ is an inverse
local time process, i.e.:
\[
\xi_l=\inf\{t;L_t>l\},
\]
where $\PR{L_t}{t\geq0}$ is (a choice of) the local time at 0 for a
\emph{radial Ornstein--Uhlenbeck process}, with dimension $\delta=2
(1-\a
)\in(0,2)$ and parameter $\mu>0$. This process is defined as the
square-root of the $\mathbb{R}_+$-valued diffusion $\PR{Z_t}{t\geq0}$
which solves:
\[
Z_t=2 \int_0^t
\sqrt{Z_s} \,\mathrm{d}\beta_s-2\mu\Io {0} {t}
{Z_s} {s}+\delta t
\]
where $\PR{\beta_s}{s\geq0}$ denotes a standard real-valued Brownian
motion starting from 0. This family of subordinators $\xi$ was studied
in Pitman--Yor \cite{PY}. We have also devoted some study to this
process in \cite{HY}.

By Pitman--Yor \cite{PY} (see also Hirsch--Yor \cite{HY}), the
Laplace--Bernstein exponent of $\xi$ is, for a suitable choice of the
local time,
\begingroup
\abovedisplayskip=6.5pt
\belowdisplayskip=6.5pt
\[
\Phi(s)=\frac{\Gamma (\frac{s}{2\mu}+\a )}{\Gamma
(\frac
{s}{2\mu} )}
\]
and
\[
a=0\quad\mbox{and}\quad\l(\mathrm{d}x)=\frac{2 \mu\a}{\Gamma
(1-\a)} \frac{\mathrm{e}^{-2\mu\a x}}{(1-\mathrm{e}^{-2\mu x})^{1+\a}}
\,\mathrm{d}x.
\]
Clearly, $\Phi$ is a complete Bernstein function, but $\Phi\notin
\mathcal{S}$. By Corollary \ref{c31}, $\mathcal{I}$ is m.i.d. A~simple
computation, from formula (\ref{e41}), yields:
\[
\kappa(\mathrm{d}x)=\frac{1-\mathrm{e}^{-2\mu\a x}}{1-\mathrm
{e}^{-2\mu
x}} \,\mathrm{d}x.
\]
Then, with the notation of Proposition \ref{p35}, one has:
\[
\ell(x)=\frac{\mathrm{e}^{-2\mu\a x}-\mathrm{e}^{-2\mu
x}}{(1-\mathrm
{e}^{-2\mu x}) (\mathrm{e}^x-1)}.
\]
It is not difficult to see that $\ell$ is a decreasing function. Then,
by Proposition \ref{p35}, $\log\mathcal{I}$ is self-decomposable.
\begin{itemize}
\item
If $2 \mu\a=1$, then by Theorem \ref{t21} and formula (\ref{e23}),
\[
R(r)=\frac{\Gamma(\a r)}{\Gamma(\a)}\quad\mbox{and}\quad I(r)=\frac
{\Gamma(\a) \Gamma(r)}{\Gamma(\a r)}.
\]
Consequently, $\mathcal{R}\st\gamma_{\a}^{\a}$ and the law of
$\mathcal
{I}$ is: $( \mathbb{E}[\tau_{\a}^{-\a}])^{-1} x \mathbb{P}_{\tau
_{\a
}^{-\a}}(\mathrm{d}x) $, where $\tau_{\a}$ denotes a standard $\a
$-stable positive variable, i.e.,
\[
\mathbb{E}\bigl[\exp(-s \tau_{\a})\bigr]=\exp\bigl(-s^{\a}
\bigr),
\]
and $\mathbb{P}_{\tau_{\a}^{-\a}}(\mathrm{d}x)$ denotes the law of
$\tau_{\a}^{-\a}$.

Besides, with the notation of Proposition \ref{p35}, one has:
\[
\jmath(x)=\frac{\mathrm{e}^{-x}}{1-\mathrm{e}^{-2\mu x}}.
\]
Clearly, $\jmath$ is a decreasing function. Then, by Proposition \ref
{p35}, $\log\mathcal{R}$ also is a self-decomposable variable.
\item
If $2 \mu(1-\a)=1$, then by Theorem \ref{t22} and formula (\ref{e23}),
\begin{eqnarray*}
&\displaystyle I(r)=(1-\a)^{1-r} \Gamma\bigl((1-\a) (r-1)+1\bigr)&
\\
&\displaystyle\mbox{and }\quad R(r)=\frac{(1-\a)^{r-1} \Gamma(r)}{\Gamma((1-\a
) (r-1)+1)}.&
\end{eqnarray*}
Consequently, $\mathcal{I}\st(1-\a)^{-1} \mathrm{\mathbf{e}}^{1-\a}$
and $\mathcal{R}\st(1-\a) \tau_{1-\a}^{\a-1}$.\vadjust{\goodbreak}

Besides, with the notation of Proposition \ref{p35}, one has:
\[
\jmath(x)=\frac{1}{\mathrm{e}^{x}-1}-\frac{1}{\mathrm{e}^{\frac
{1}{1-\a}x}-1}.
\]
It is not difficult to show that $\jmath$ is a decreasing function.
Then, by Proposition \ref{p35}, $\log\mathcal{R}$ also is a
self-decomposable variable.

Note that, in this case, the Bernstein function $\Phi$ writes:
\[
\Phi(s)=\frac{\Gamma((1-\a) (s-1)+1)}{\Gamma((1-\a) s)},
\]
and therefore, $\Phi$ is the Bernstein function of example \ref{sss432}
with $c=\a^{-1}$, after changing $\a$ into $1-\a$.
\end{itemize}
\endgroup
%

\section{Relating our results to Urbanik's}\label{s5}
As we wrote in the Introduction, some of our main results, notably
those of Section \ref{s3}, may also be found in Urbanik's paper \cite
{U} (for convenience, we shall simply write \U$ $ when referring to
this paper). However, in order that our reader may have some cosiness
in comparing our results with those in \U, we first need to recall and
explain the main notation in \U, which we undertake in the next
Subsection \ref{ss51}, while Subsection \ref{ss52} is devoted to the
statement and explanation (both in \U's manner, and with our notation)
of the main relevant results in \U. Finally, in Subsection~\ref{ss53},
we compare some of \U's results with ours.
\subsection{Basic notation in \U}\label{ss51}
\begin{enumerate}[a)]
\item[a)] In \U, a subordinator $\PR{\xi_l}{l\geq0}$\footnote{We keep
using, whenever convenient, our notation and we confront it / compare
it / with that in \U.}, and any quantity related to it, is always
referred to, or tagged, by its representing measure $M(\mathrm{d}x)$ on
$\R_+$, as follows:
\[
\mathbb{E}\bigl[\exp(-z \xi_l)\bigr]=\exp\bigl(-l \Phi(z)\bigr) ,
\quad l,z\geq0,
\]
where
%
\begin{equation}
\label{e51} \Phi(z)=\I{\R_+} {\frac{1-\mathrm{e}^{-zx}}{1-\mathrm{e}^{-x}}} {M} {x},
\end{equation}
with the function: $ {x\longrightarrow\frac{1-\mathrm
{e}^{-zx}}{1-\mathrm{e}^{-x}}}$ being taken equal to $z$ for $x=0$. (In
fact, in \U, the Bernstein function $\Phi(z)$ is denoted: $\langle M\rangle (z)$.)
Note that (\ref{e51}) writes:
%
\begin{equation}
\label{e52} \Phi(z)=z M\bigl(\{0\}\bigr)+\I{(0,+\infty)} {\frac{1-\mathrm
{e}^{-zx}}{1-\mathrm
{e}^{-x}}} {M}
{x}.
\end{equation}
Clearly, from (\ref{e52}), the L\'evy measure -- which we shall denote
as $\l(\mathrm{d}x)$ -- on $(0,+\infty)$, associated to $\PR{\xi_l}{l\geq0}$ is
%
\begin{equation}
\label{e53} \l(\mathrm{d}x)=1_{(0,+\infty)}(x) \bigl(1-\mathrm{e}^{-x}
\bigr)^{-1} M(\mathrm{d}x).
\end{equation}
However, in \U, the L\'evy measure is never mentioned, all reference
being made to $M(\mathrm{d}x)$, which is a \emph{bounded} measure on
$\R_+$.
\item[b)] As an important example of the notation in \U$ $ which we
presented so far, we note (\U, p.~494) that the representing measure
of the standard gamma process $\PR{\gamma_l}{l\geq0}$ is denoted as:
%
\begin{equation}
\label{e54} \Pi(\mathrm{d}x)=\frac{\mathrm{e}^{-x} (1-\mathrm{e}^{-x})}{x} \,\mathrm{d}x.
\end{equation}
Indeed, the L\'evy measure $\l_{\gamma}(\mathrm{d}x)$ of the gamma
process $\PR{\gamma_l}{l\geq0}$ is:
\[
\l_{\gamma}(\mathrm{d}x)=\frac{\mathrm{e}^{-x}}{x} \,\mathrm{d}x.
\]
Also, formulae (1.3) and (1.4) in \U$ $ may read:
%
\begin{equation}
\label{e55} \Phi_{\gamma}(z)\equiv\langle \Pi\rangle (z)=\log(1+z).
\end{equation}
\item[c)] In order to understand the main results in \U, and their
connection with ours, we still need to recall concepts related to the
$S$-transform on positive bounded measures on $\R_+$:
$M\longrightarrow
SM$ as introduced in \U. Rather than describing the $S$-transform
directly, we prefer first to indicate how it acts on Bernstein
functions. In \U, formula (1.21), we find:
%
\begin{equation}
\label{e56} \Phi_{SM}(z)=\log \bigl(\Phi_{\overline{M}}(1+z) \bigr)
\end{equation}
where: $\overline{M}=(M(\R_+))^{-1} M$. It is easily seen, and this is
confirmed by \U, formula (1.20), that the subordinator $\PR{\xi_{SM}(l)}{l\geq0}$ may be constructed as $\PR{\widetilde{\xi
}_{\gamma
_{l}}}{l\geq0}$, that is by $\gamma$-subordination of the subordinator
$\PR{\widetilde{\xi}_l}{l\geq0}$, which itself is the Esscher
transform of the subordinator $\xi_{\overline{M}}$; precisely,
\[
\mathbb{P}^{\widetilde{\xi}}_{ | \mathcal{F}_l}=\exp (-X_l+l)\cdot
\mathbb{P}^{\xi_{\overline{M}}}_{  |\mathcal
{F}_l}
\]
where $\mathbb{P}^{\widetilde{\xi}}$ and $\mathbb{P}^{\xi
_{\overline
{M}}}$ are the laws of the subordinators on the canonical Skorokhod
space $D([0,+\infty))$, where $X_l(\omega)=\omega(l)$, and $\mathcal
{F}_l=\sigma\{X_m;m\leq l\}$. Finally, for concreteness, let us give
the explicit form of the measure $SM$; this is (\U, formula (1.25)):
%
\begin{equation}
\label{e57} SM(\mathrm{d}x)=\bigl(1-\mathrm{e}^{-x}\bigr)
\mathrm{e}^{-x} \Io{0} {\infty } {l^{-1} e_+(l M) (
\mathrm{d}x)} {l},
\end{equation}
where $e_+(l M)$ is the notation in \U$ $ for the law of $\xi_{M}(l)$, for given $l$ (which we denoted as $\mu_l$ in Subsection
\ref{ss21}).
\end{enumerate}
%
\begin{rem}\label{r51}
A simple computation from (\ref{e57}) shows that, for $s> 0$,
\[
\I{\R_+} {\mathrm{e}^{-sx} \bigl(1-\mathrm{e}^{-x}
\bigr)^{-1} \mathrm{e}^{x} x} {SM} {x}=\frac{\Phi_M'(s)}{\Phi_M(s)}.
\]
Hence, the measure: $(1-\mathrm{e}^{-x})^{-1} \mathrm{e}^{x} x
SM(\mathrm{d}x)$ is our measure $\kappa$ associated to $\Phi_M$ (see
Proposition~\ref{p21}).
\end{rem}
%

\subsection{Main results in \U}\label{ss52}
We are now in a right position to state clearly and precisely the main
results in \U, before comparing them to ours.
\begin{enumerate}[a)]
\item[a)] We first write:
\[
\frac{1}{1+z}=\frac{1}{\Phi_{\overline{M}}(1+z)} \biggl(\frac{\Phi_{\overline{M}}(1+z)}{1+z} \biggr),
\]
and note that the left-hand side, which is the Laplace transform of a
standard exponential variable $\mathrm{\mathbf{e}}$, appears on the
right-hand side as a product of two Laplace transforms. Precisely,
there is the following
theorem.
%
\begin{thmm}[(\U, Proposition 2.1, Theorem 2.3)]\label{t51}
\textup{i)} To every positive, bounded measure~$M$, we may associate a
positive r.v. $Y_M$ such that:
\begin{enumerate}[ii)]
%
%
\begin{equation}
\label{e58} \mathrm{\mathbf{e}}\st\xi_{SM}(1)+Y_M,
\end{equation}
where, in the right-hand side $Y_M$ is independent of $\xi_{SM}(1)$.
\item[ii)]$Y_M$ is infinitely divisible if and only if $SM\leq\Pi$, in
which case, there is the identity in law:
\[
Y_M\st\xi_{\Pi-SM}(1)
\]
and, of course the identity in law (\ref{e58}) extends at the level of
subordinators:
%
\begin{equation}
\label{e59} \PR{\gamma_l} {l\geq0}\st\PR{\xi_{SM}(l)+
\xi_{\Pi-SM}(l)} {l\geq0}
\end{equation}
where, in the right-hand side, the two subordinators are independent.
\end{enumerate}
\end{thmm}
Besides, in \U, Theorem 2.7, a characterization of the measures $\Pi
-SM$, assumed to be $\geq0$, is given.
\item[b)] We are now -- almost -- in the right position to state and prove
another important result found in \U, Theorem 3.3, pertaining to the
m.i.d. property of $\mathcal{I}_M$.
%
\begin{thmm}[(\U, Theorem 3.1 and Theorem 3.3)]\label{t52}
\textup{i)} Set:
\[
\mathcal{I}_M=\Io{0} {\infty} {\exp\bigl(-\xi_M(l)
\bigr)} {l}.
\]
There exists a random variable $\mathcal{R}_M$, independent of
$\mathcal
{I}_M$, such that:
\[
\mathrm{\mathbf{e}}\st\mathcal{I}_M\cdot\mathcal{R}_M.
\]
\begin{enumerate}[ii)]
\item[ii)]$\mathcal{R}_M$ is m.i.d. and
\[
-\log\mathcal{R}_M\st e(r_M,SM) \quad\mbox{with}\quad
r_M=-\log M(\R_+)+l(SM).
\]
\item[iii)]$\mathcal{I}_M$ is m.i.d. if and only if: $SM\leq\Pi$.
\item[iv)] When the condition: $SM\leq\Pi$ is satisfied, then:
\[
-\log\mathcal{I}_M\st e(i_M,\Pi-SM) \quad\mbox{with}
\quad i_M=\log M(\R_+)+l(\Pi-SM).
\]
\end{enumerate}
\end{thmm}
In order that the statement of Theorem \ref{t52} be totally
understandable for our reader, it only remains to explain about the
notation $e(a,M)$ and $l(M)$, which is featured twice (in~ii) and~iv))
in the above Theorem \ref{t52}.

The notation $e(a,M)$ indicates a real-valued r.v., or rather its law,
which is infinitely divisible, and whose characteristic function
$\widetilde{e}(a,M)$, in the notation of \U, is given by:
\[
\widetilde{e}(a,M) (s)=\exp \biggl[\mathrm{i}as+\int_{\R_+}
\biggl(\mathrm {e}^{\mathrm{i}sx}-1-\frac{\mathrm{i}sx}{1+x^2} \biggr)\frac
{M(\mathrm
{d}x)}{(1-\mathrm{e}^{-x})^2}
\biggr],
\]
whereas $l(M)$ is a real number defined as:
\[
l(M)=\int_{\R_+} \biggl(\mathrm{e}^{-x}-1+
\frac{x}{1+x^2} \biggr)\frac
{M(\mathrm{d}x)}{(1-\mathrm{e}^{-x})^2}.
\]
\end{enumerate}
%
\subsection{Comparison of \U's results with ours}\label{ss53}
We are now able to see that our main results of Section \ref{s3} appear
in the above Theorem~\ref{t52}. Indeed, by Remark \ref{r51},
%
\begin{equation}
\label{e510} SM(\mathrm{d}x)=\frac{1-\mathrm{e}^{-x}}{x \mathrm{e}^x} \kappa (\mathrm{d}x)
\end{equation}
where $\kappa$ still denotes the measure whose Laplace transform is
$\Phi'_M/\Phi_M$. In particular, by formula~(\ref{e54}), \U's
condition: $SM\leq\Pi$ is equivalent to our condition: $\kappa
(\mathrm
{d}x)\leq\mathrm{d}x$. On the other hand, we also obtain from (\ref{e510}):
%
\begin{equation}
\label{e511} \widetilde{e}\bigl(l(SM),SM\bigr) (s)=\exp \biggl[\int
_{(0,+\infty)}\frac
{\mathrm
{e}^{\mathrm{i}sx}-1+\mathrm{i} s (\mathrm{e}^{-x}-1)}{x (\mathrm
{e}^{-x}-1)} \kappa(\mathrm{d}x) \biggr].
\end{equation}
Moreover,
%
\begin{equation}
\label{e512} \widetilde{e}\bigl(l(\Pi),\Pi\bigr) (s)=\exp \biggl[\int
_0^{+\infty}\frac
{\mathrm
{e}^{\mathrm{i}sx}-1+\mathrm{i} s (\mathrm{e}^{-x}-1)}{x (\mathrm
{e}^{-x}-1)} \,\mathrm{d}x \biggr]
\end{equation}
and, by (\ref{e52}),
%
\begin{equation}
\label{e513} M(\R_+)=\Phi_M(1).
\end{equation}
Thus, in view of the above formulae (\ref{e511}), (\ref{e512}) and
(\ref
{e513}), property ii) in Theorem \ref{t52} corresponds to our Theorem
\ref{t31}, while properties iii) and iv) in Theorem \ref{t52}
correspond to our formula (\ref{e35}), Theorem \ref{t32} and Theorem
\ref{t33}.

Moreover, in \U, Theorem 2.5, it is stated that, if $\Phi_M$ is a
complete Bernstein function, then $\mathcal{I}_M$ is m.i.d., which is
our Corollary \ref{c31}.
%

%


\printhistory

\end{document}